\documentclass[twocolumn]{autart}    

\usepackage{graphics} 
\usepackage{epsfig} 
\usepackage{amsmath} 
\usepackage{amssymb}  
\usepackage{breqn}
\usepackage{algorithm}
\usepackage{algpseudocode}

\newcommand{\algrule}[1][.2pt]{\par\vskip.5\baselineskip\hrule height #1\par\vskip.5\baselineskip}

\usepackage{booktabs}
\usepackage{subcaption}
\usepackage{caption}

\newtheorem{definition}{Definition}
\newtheorem{proposition}{Proposition}

\begin{document}

\begin{frontmatter}

\title{On the Design of Rational Polynomial State Feedback Controllers\thanksref{footnoteinfo}} 

\thanks[footnoteinfo]{Matthew Newton and Zuxun Xiong contributed equally. This work was supported by EPSRC grants EP/L015897/1 (to M. Newton) and EP/Y014073/1 \& EP/X031470/1 \& UKRI2108 (to A. Papachristodoulou) and the Worcester College Oxford Tony Corner Research Fund.}

\author[UK]{Matthew Newton}\ead{matthew.ant.newton@gmail.com},    
\author[UK]{Zuxun Xiong}\ead{zuxun.xiong@eng.ox.ac.uk},               
\author[Switzerland]{Han Wang}\ead{: hanwang@control.ee.ethz.ch},
\author[UK]{Antonis Papachristodoulou}\ead{: antonis@eng.ox.ac.uk}

\address[UK]{Department of Engineering Science, University of Oxford, OX1 3PJ, Oxford, UK}  
\address[Switzerland]{Automatic Control Laboratory (IfA), ETH Z\"urich, Z\"urich, Switzerland}             


\begin{abstract}                          
One of the desirable objectives in feedback control design is to formulate and solve the design problem as an optimisation problem that is convex, so that an optimal solution can be found efficiently. Unfortunately many control design problems are non-convex: approximations, relaxations, or iterative schemes are usually employed to solve them. Several such approaches have been developed in the literature, for example Sum-of-Squares (SOSs) methods have been used for systems described by polynomial dynamics. Alternatively, and relevant to this paper, one can choose a (non-unique) linear-like representation of the system and solve the resulting state-dependent Linear Matrix Inequalities (LMIs) or use SOSs optimisation techniques to derive a control law. This SOS method has been shown to effectively design polynomial and rational controllers for nonlinear polynomial systems, offering a broader class of controllers and the potential for improved performance and robustness guarantees. In this paper, we start off by considering rational functions as controllers for nonlinear systems and propose a procedure for designing such controllers by iteratively solving convex SOS optimisation problems. Our approach decouples the controller structure from the system dynamics and incorporates it as a constraint within the optimisation problem, which results in an optimisation that co-designs aspects of the controller and the Lyapunov function at the same iterative step. We theoretically establish the properties of this procedure by showing that several existing rational controller design methods can be recovered as special cases of this procedure. The proposed method is evaluated on various nonlinear benchmark system examples, demonstrating improved performance and robustness over both polynomial controllers and rational controllers obtained by existing approaches.
\end{abstract}

\end{frontmatter}

\section{Introduction}\label{sec:rationalintro}
Synthesising state-feedback controllers is a key task in control systems design. In the case of linear systems, linear state feedback controllers have been widely studied due to their simplicity and the availability of convex optimisation-based design methods. There is a comprehensive literature on how the resulting Lyapunov conditions can be transformed into Linear Matrix Inequalities (LMIs), thereby formulating linear controller synthesis as a convex optimisation problem~\cite{boyd1994linear}. The more general problem of co-designing a Lyapunov function and a nonlinear control law for a nonlinear system is more challenging; and linear controllers exhibit limitations when applied to systems with significant nonlinearities, constraints and operational mode switching~\cite{nikolaou2003linear}. Various methods have therefore been developed to address this issue, most of which fall into two main categories: local approximations and iterative procedures. For instance, linearisation may be employed to construct convex approximations around an operating point~\cite{khalil2002nonlinear}, but this makes the solutions only valid locally. Iterative schemes have also been proposed to solve non-convex optimisation problems by a series of convex approximations particularly when the formulation involves special structures such as bilinear terms~\cite{doelman2016sequential}. However, such methods usually have no optimality guarantees, and may be sensitive to initial conditions. This has motivated the development of more advanced iterative frameworks, such as sequential convex programming, which provide local convergence guarantees for the non-convex problems arising in nonlinear system analysis~\cite{cunis2022sequential}. 

The development of convex formulations for nonlinear control analysis problems - notably the Sum of Squares (SOSs) approach - has allowed researchers to reconsider the synthesis question. The SOS hierarchy provides a tractable relaxation method to formulate NP-hard polynomial positivity problems into tractable convex problems~\cite{parrilo2000structured}. This method has enabled a wide range of analysis tasks, such as constructing Lyapunov functions for nonlinear stability analysis~\cite{apap1} and barrier functions for safety verification~\cite{prajna2007framework}. Building on this, the framework has been used to verify closed-loop properties of nonlinear systems such as stability~\cite{korda2017stability} and inner-approximations of the region of attraction (ROA)~\cite{topcu2008local,korda2014controller}. Moreover, the SOS approach opened the door to previously difficult nonlinear synthesis tasks which can be solved efficiently using convex optimisation - early works demonstrated its potential for a variety of control design applications~\cite{jarvis2003some}. This led to the development of structured methodologies, such as leveraging dissipation inequalities to formulate convex design problems~\cite{ebenbauer2006analysis}, and exploiting specific system structures to decompose the non-convex synthesis problem into a tractable convex problem~\cite{prajna2004nonlinear2,prajna2004nonlinear}.

The SOS technique has enabled tractable synthesis of nonlinear state feedback controllers, but these were primarily polynomial functions of the state. However, such controllers often perform poorly far from the operating region, producing excessively large control inputs when the state deviates from equilibrium, potentially leading to instability or actuator saturation which can introduce discontinuous cut-offs. Alternatively, rational controllers have the potential to offer a more expressive and structurally flexible alternative. First, the structural form of rational functions provides a principled means to address actuator constraints, enabling more effective handling of input saturation. Moreover, rational functions exhibit stronger approximation capabilities than polynomials of comparable complexity. Recent studies have shown that rational functions can approximate complex nonlinear mappings more accurately and efficiently than neural networks (NN) with the same number of parameters~\cite{mtel,vpei}. This combination of input regularisation and expressive power makes rational controllers particularly appealing in addressing nonlinear control problems. For instance, in aerial vehicle trajectory tracking tasks, rational controllers show superior performance compared to other traditional control strategies~\cite{ktan2}.

Motivated by these insights, this work develops a new iterative procedure for synthesising rational controllers, leveraging their structural advantages and strong approximation properties while ensuring a tractable synthesis process through convex optimisation.

\subsection{Related Work}
Rational controllers have the potential to address difficult tasks, but their design is more complicated; challenges include issues such as convexity, restrictions on the dynamics, and flexibility in the co-design of a rational controller and a Lyapunov function. 

Rational controllers designed within the SOS framework have shown effectiveness when applied to bilinear systems~\cite{vatani2014control,xie2025bilinear,strasser2024koopman}. In particular, \cite{vatani2014control} shows that for an unstable discrete-time bilinear system, only rational controllers can achieve global stabilisation. However, the controller design in that work is based on a pre-defined Lyapunov function, which limits design flexibility. In~\cite{xie2025bilinear} and~\cite{strasser2024koopman}, SOS techniques were applied to design stable rational controllers under a robust MPC framework and a Koopman bilinearisation framework, respectively. A common drawback of both approaches is that the denominator of the rational controller is fixed \emph{a priori}, allowing only the numerator to be optimised. 

For more general control-affine systems, the well-known Sontag's formula provides a universal construction for stabilising controllers based on analytical cancellation techniques~\cite{eson}. Although the resulting controller is an algebraic function of the state, rather than a strictly rational polynomial one, its structure has inspired subsequent research in rational controller design. This formula can also be easily extended and applied to cases with bounded input constraints~\cite{ylin}. However, this method assumes a known Lyapunov function. To overcome this, some works consider the joint synthesis of the controller and the Lyapunov certificate, formulating the resulting problem using the SOS framework~\cite{prajna2004nonlinear,mvat,pylorof2016analysis,huang2013robust,madeira2025nonlinear}. 

To address the non-convexity introduced by joint synthesis, a widely used approach is the iterative method, where the Lyapunov function and the controller are designed alternately~\cite{mvat,pylorof2016analysis}. While this strategy improves tractability, it often sacrifices flexibility in the resulting controller, as will be discussed in detail in later sections. Few notable works solve this problem to some extent~\cite{prajna2004nonlinear,huang2013robust,madeira2025nonlinear}. In~\cite{prajna2004nonlinear,huang2013robust}, the rational feedback law emerges from a specific parametrisation of the inverse Lyapunov matrix, but the resulting controller shares a fixed denominator with the Lyapunov function. This tight coupling limits design flexibility. The method proposed in~\cite{madeira2025nonlinear} adopts a dissipativity-based framework for static output-feedback design, where the rational structure arises from a matrix inverse in the stability condition. This condition introduces an additional SOS constraint that must be enforced at each iteration. Furthermore, the denominator is restricted to a product of diagonal terms shared across input channels, limiting flexibility in the controller structure. An alternative line of research bypasses the challenges of joint synthesis by leveraging a different stability criterion based on density functions, which transforms the synthesis problem into a convex SOS program from which a rational controller can be designed~\cite{prajna2004nonlinear2}. Alternatively, \cite{korda2016controller} proposes a hierarchy of SOS programs based on occupation measures which yields a sequence of asymptotically optimal rational controllers.

\subsection{Our Contribution}
To address the existing limitations, most notably the restriction to control-affine systems, the non-convexity associated with joint synthesis of the controller and Lyapunov function, and the structural limitation on the controller form, we look at this synthesis problem from a different perspective and propose a new procedure for the design of rational controllers based on SOS. The contributions of this work are as follows:
\begin{itemize}
    \item We reformulate the synthesis problem by decoupling the controller expression, embedding it as a constraint and setting the controller input as an independent variable in the optimisation problem. These steps enable structural decoupling between the controller and Lyapunov function, making the problem amenable to convex optimisation via SOS techniques.
    \item Building on the above formulation, we develop an iterative design procedure that alternates between controller parameters and multiplier updates. Each iteration is a convex optimisation problem; importantly it enables (in every other step) the joint optimisation of the rational controller and the Lyapunov function within the same procedure, which can enlarge the admissible solution space and relax some of the structural restrictions present in traditional iterative methods. 
    \item We incorporate performance-oriented criteria such as the Lyapunov decay rate and the size of the Region of Attraction (ROA) into the design procedure, thereby enhancing its versatility and enabling performance-aware controller synthesis.
    \item Through numerical examples, we demonstrate that the proposed SOS-based method can effectively synthesise stable rational controllers for a variety of unstable nonlinear systems with good performance.
\end{itemize}





\section{Preliminaries} \label{sec:preliminary}

\subsection{Rational Functions} 
\begin{figure}[!ht] 
    \centering
\includegraphics[width=0.3\textwidth]{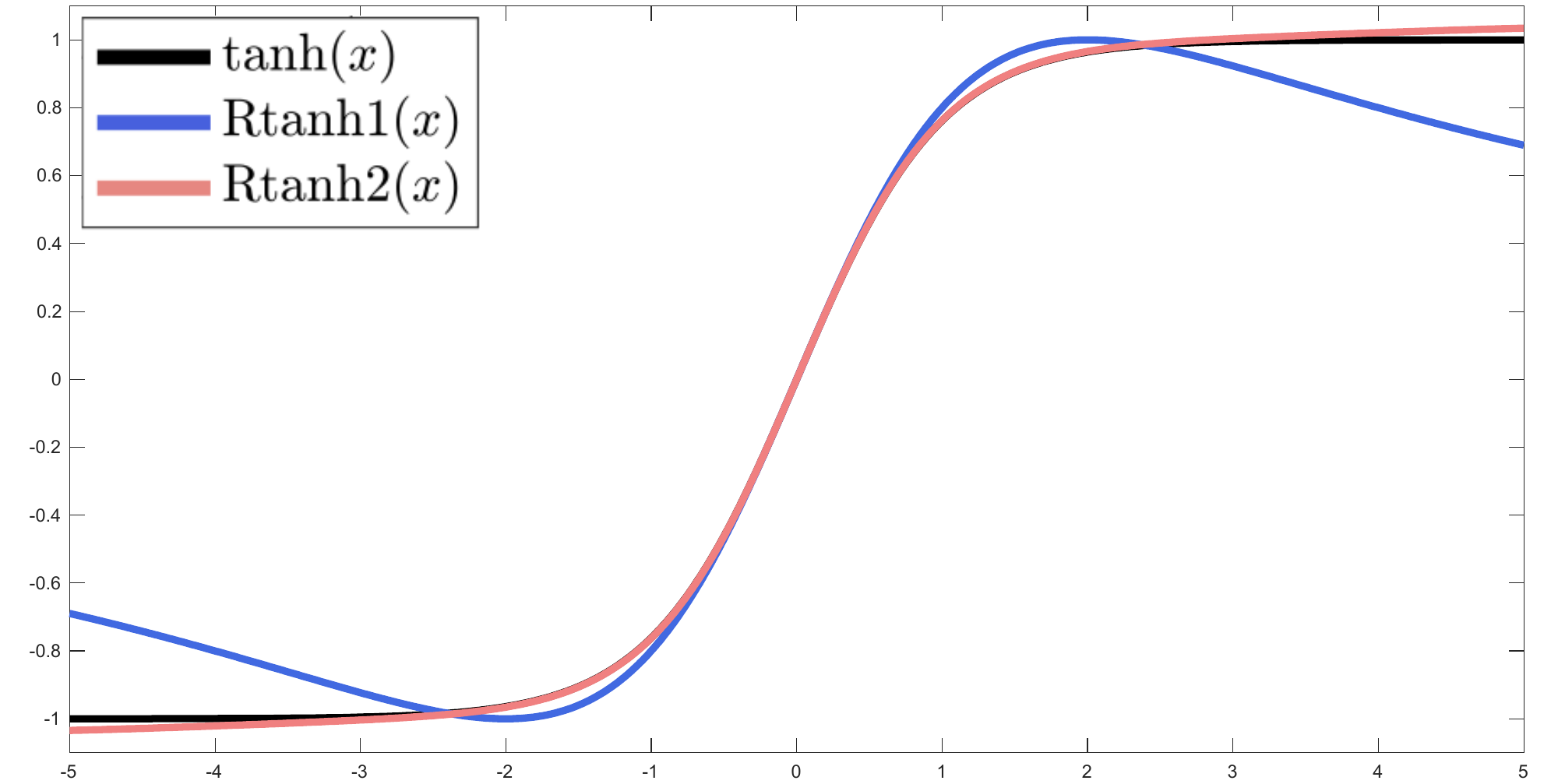}
        \caption{Comparison between $\mathrm{Rtanh}_1(x)$, $\mathrm{Rtanh}_2(x)$ and $\mathrm{tanh}(x)$}
        \label{Ex-tanh-error}
\end{figure}

Here we consider rational functions in the form $\frac{p(x)}{q(x)}$, where $p(x)$ and $q(x)$ are polynomials. Rational functions can be used to approximate a wide class of other nonlinear functions. For example, the simple rational expression
$\mathrm{Rtanh}_{1}(x) = \frac{4x}{x^{2} + 4}$, is a good approximation of the nonlinear function $\mathrm{tanh}(x) = \frac{e^{x} - e^{-x}}{e^{x} + e^{-x}}$, close to the origin. To further improve the accuracy, we apply the AAA method \cite{nakatsukasa2018aaa} to obtain a fourth degree rational function $\mathrm{Rtanh}_{2}(x) = \frac{1.7x^3+22.7x}{0.1x^4+9.2x^2+22.7}$, which can approximate $\mathrm{tanh}(x)$ quite well. These functions are shown in Fig. \ref{Ex-tanh-error}: $\mathrm{tanh}(x)$ is almost fully covered by $\mathrm{Rtanh}_{2}(x)$.

\begin{figure}[!ht] 
    \centering  
    \includegraphics[width=0.8\linewidth]{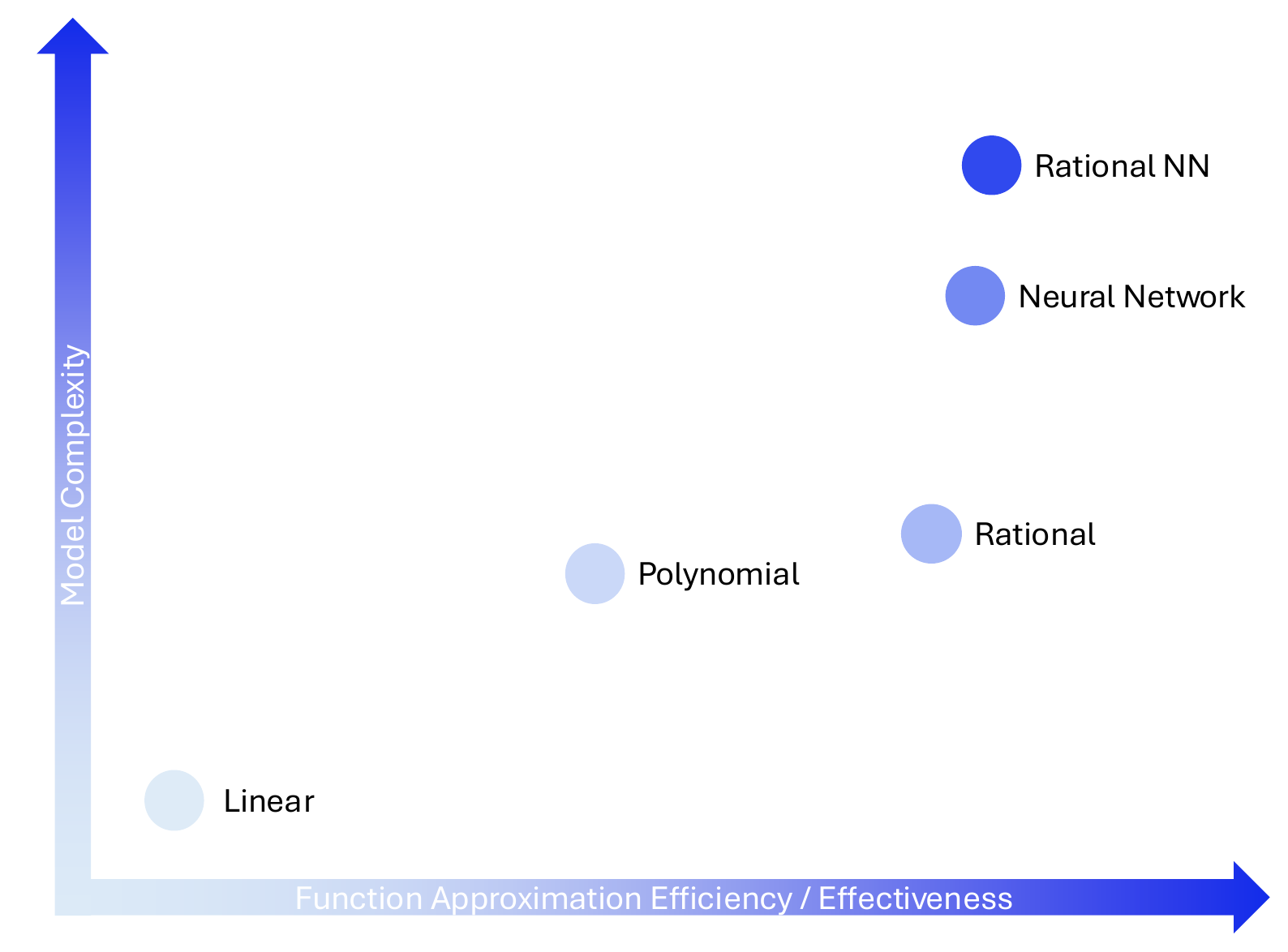}
    \caption{Function Classes: Approximation Power vs. Model Complexity.} \label{complexity-appro}
\end{figure}

At the same time, Neural network controllers have recently been introduced in the control toolbox, given that Neural Networks are regarded as powerful universal function approximators due to their nonlinear activation functions and large number of internal parameters~\cite{hornik1989multilayer}. However, this expressive power comes at the cost of increased model complexity and difficulties in verification and control design. Rational functions are a middle ground between low complexity functions and high complexity NNs and are therefore a great candidate for feedback controllers. Fig.~\ref{complexity-appro} illustrates a qualitative comparison of representative function classes in terms of their approximation power and overall model complexity. The horizontal axis reflects the expressive capacity or efficiency in approximating general functions, while the vertical axis uses `complexity' to loosely capture aspects such as the number of parameters, architectural or implementation complexity, and computational or design effort. As shown, rational functions offer a favourable trade-off, achieving strong approximation capabilities with moderate complexity. Moreover, rational NNs -- which use rational functions to replace standard activations -- can further extend expressivity while preserving some of the structural advantages of rational approximators~\cite{nbou,mnew6}. We emphasise that this diagram is intended as a conceptual tool rather than a quantitative chart; the actual complexity and approximation efficiency may vary across different cases and implementations. For more systematic comparisons, the reader is referred to~\cite{mtel,vpei,alimov2021geometric}.

Another attractive property of rational controllers over polynomial ones lies in their inherent ability to regulate and even bound the control magnitude at large state values through their denominator structure, thereby achieving a continuous closed-loop response without relying on discontinuous saturation functions.



\subsection{Sum of Squares Programming} \label{sec:sosprogpsatz}
We first review some background that is relevant for the rest of the paper. Denote by $\mathbb{R}[x_{1}, \dots , x_{n_x}]$ the set of polynomials in $x_1, \ldots, x_{n_x}$ with real coefficients; we will use $x = (x_1, \ldots, x_{n_x})$ for simplicity.
\begin{definition} \label{def:sos}
    A polynomial $p(x) \in \mathbb{R}[x]$ is said to be a \emph{Sum of Squares} if it can be expressed as 
    \begin{equation*}
        p(x) = \sum_{i=1}r_{i}^{2}(x), 
    \end{equation*}
    where $r_i(x) \in \mathbb{R}[x]$.
\end{definition}
 We denote the set of polynomials that admit this decomposition by $\Sigma[x]$ and say `$p(x)$ is SOS'. It is obvious that if a polynomial $p(x)$ admits an SOS decomposition, it is non-negative for all $x \in \mathbb{R}^{n_x}$ -- the converse is not true, and in fact verifying if $p(x)$ is non-negative is NP-hard, while verifying if $p(x)$ is SOS can be formulated as a Semidefinite Program.  
This can be achieved with toolboxes such as SOSTOOLS~\cite{sostools}, SOSOPT~\cite{seiler2013sosopt}, and GloptiPoly~\cite{henrion2003gloptipoly} in MATLAB. For additional details on SOS programming the reader is referred to~\cite{apap1,jarvis2003some,chesi2011domain,anderson2015advances}. 


\subsection{Lyapunov Stability Analysis}
In this section we outline how SOS programming can be used to verify the stability of a system under constraints. Consider a continuous-time system
\begin{equation} \label{eq:statespace1}
\begin{aligned}
    \dot{x}(t) &= f(x(t),u(t)), \\
    g_{i}(x(t),u(t)) &\geq 0, \: \forall i = 1, \dots, q_{1}, \\
    h_{j}(x(t),u(t)) &= 0, \: \forall j = 1, \dots, q_{2},
\end{aligned}
\end{equation}
where $x(t) \in \mathbb{R}^{n_{x}}$ and $u(t) \in \mathbb{R}^{n_{u}}$ are the system states and inputs respectively, $f(x(t),u(t))$ is the plant model, which is a vector of polynomial or rational functions in $(x(t),u(t))$ without singularity in a set $S$: a semi-algebraic set defined as $S:=\{(x(t),u(t)) \mid g_i(x(t),u(t))\geq 0, h_j(x(t),u(t))=0, \forall i=1, \ldots, q_1, j=1, \ldots, q_2\}$, where $g_{i}(x(t),u(t)) \geq 0$ and $h_{j}(x(t),u(t)) = 0$ are polynomial constraints on the system. Throughout this paper we drop the time dependence on the system states to improve readability. Based on a central theorem of real algebraic geometry, Positivstellensatz~\cite{gsten1}, which provides a relationship between an algebraic condition and the emptiness of a semi-algebraic set, a Lyapunov stability SOS condition is provided in following proposition:

\begin{proposition}[\protect {\cite[Proposition 1]{mnew5}}]\label{prop:ss1stability}
Consider System \eqref{eq:statespace1} and suppose there exists a polynomial function $V(x)$ that satisfies the following conditions
\begin{equation} \label{sosopt1}
\begin{aligned}
    &V(x) - \rho(x) \in \Sigma[x],\quad\rho(x) > 0, \\
    &-\frac{\partial V}{\partial x}(x) f(x,u) - \sum_{j=1}^{q_2}t_{j}(x,u)h_{j}(x,u) \\
    &-\sum_{i=1}^{q_1}s_{i}(x,u)g_{i}(x,u) 
    \in \Sigma[x,u], \\
    &s_{i}(x,u) \in \Sigma[x,u], \: \forall i = 1, \dots, q_{1}, \\
    &t_{j}(x,u) \in \mathbb{R}[x,u], \: \forall j = 1, \dots, q_{2}, 
\end{aligned}
\end{equation}
then the equilibrium of the system \eqref{eq:statespace1} is stable.
\end{proposition}

\section{Design Paradigms for Rational Nonlinear Controllers}\label{sec:traditional iterative method}
In this section, we summarise two major approaches to the design of rational controllers for nonlinear systems. The first is based on analytical cancellation using a known Lyapunov function, inspired by Sontag's universal formula. The second class of methods aims to jointly synthesise a Lyapunov function and a rational controller by solving a nonconvex optimisation problem, for which iterative procedures are commonly adopted. As mentioned previously, these methods primarily target control-affine systems, whose dynamics can be written as:
\begin{align} \label{eq:statespace6}
\begin{split}
    \dot{x} = f(x) + g(x)u,
\end{split}
\end{align}
where $x \in \mathbb{R}^{n_{x}}$, $u \in \mathbb{R}^{n_u}$ and $f(0) = 0$. The Lyapunov derivative condition for this system takes the form
\begin{equation} \label{eq:lyapsimple}
    -\frac{\partial V}{\partial x}f(x) -\frac{\partial V}{\partial x}g(x)u \in \Sigma[x].
\end{equation}
\subsection{Lyapunov-Based Stable Rational Controllers Design via Analytical Cancellation}
A fundamental paradigm for designing stabilising controllers is to leverage a Control Lyapunov Function (CLF). The theoretical groundwork established by Artstein guarantees that the existence of a smooth CLF implies the existence of a stabilising controller~\cite{artstein1983stabilization}. Building upon this, Sontag provided a milestone `universal formula' which, under the additional small control property of the CLF, offers an explicit and constructive method to derive a stabilising controller that is continuous at the origin~\cite{eson}. For simplicity, here we present the results for single-input systems. Specifically, given a CLF $V(x)$,
\begin{equation*}
u = \begin{cases}
    - \dfrac{a(x) + \sqrt{a(x)^2 + b(x)^4}}{b(x)}, & \text{if } b(x) \neq 0, \\
    0, & \text{if } b(x) = 0,
\end{cases}
\end{equation*}
where $a(x) := \frac{\partial V}{\partial x}(x)f(x)$ and $b(x) := \frac{\partial V}{\partial x}(x)g(x)$, is a stabilising controller of the system. This formula is an algebraic function of the state instead of a strictly rational polynomial one due to a square root term. While~\cite{eson} also proved the existence of strictly rational controllers under certain conditions, e.g., the existence of a rational CLF, a general method for their construction was not provided, leaving this an open problem.

Inspired by the analytical cancellation principle inherent in~\cite{eson}, a controller with a simpler, strictly rational polynomial structure as shown in \eqref{eq:sontag} can be designed to satisfy the stability condition \eqref{eq:lyapsimple}.
\begin{equation}\label{eq:sontag}
    u = -\frac{\frac{\partial V}{\partial x}(x)f(x) + W(x)}{\frac{\partial V}{\partial x}(x)g(x)},
\end{equation}
where $W(x)$ is a positive-definite polynomial. This form is appealing as it directly enforces stability by cancelling terms in the Lyapunov derivative. However, this simplified structure lacks the sophisticated design of Sontag's algebraic formula thereby does not inherently guarantee continuity at the origin. This potential for discontinuities, coupled with the prerequisite of a known Lyapunov function, limits its practical application and motivates the need for more advanced, systematic design frameworks that can jointly synthesise a controller and a CLF.


\subsection{Iterative Methods for Non-Convex Design}
Another method considers the joint design of a controller and the Lyapunov function. We still use the single-input case for the purpose of illustration. If we consider a rational controller $\frac{p(x)}{q(x)}$ with $q(x) > 0$ that needs to be designed, this gives the Lyapunov stability condition
\begin{equation}\label{eq:traditional sta con}
    -\frac{\partial V}{\partial x}(x)f(x)q(x) -\frac{\partial V}{\partial x}(x)g(x)p(x) \in \Sigma[x].
\end{equation}
As the Lyapunov function $V(x)$ is multiplied with the controller $p(x)$ and $q(x)$, the condition is non-convex. A commonly-used approach to solve this problem is an iterative method, where two optimisation problems are repeatedly solved. First, the controller parameters are fixed and a Lyapunov function is found, and then the Lyapunov function is fixed and the controller parameters are optimised. Such an iterative procedure has become the prevailing strategy in recent literature, with~\cite{mvat,pylorof2016analysis} being representative examples. However, the main issue with these approaches is that the Lyapunov function and controller are designed at different steps, which can lead to poor performance. 

In a later section, we will revisit both these iterative methods and previous analytical-cancellation-based methods, showing that the proposed convex iterative procedure can recover these conventional methods as special cases.

\section{Convex Iterative Synthesis of Rational Controllers}
In this section we incorporate the ideas from Section \ref{sec:traditional iterative method} to propose a different procedure to design rational controllers and Lyapunov functions.

\subsection{SOSs Stability Conditions for Rational Controllers}
The first step is to rewrite the system and decouple the equation for the controller from the plant model. Consider system \eqref{eq:statespace1} and a rational controller $\pi(x)=[\frac{p_1(x)}{q_1(x)},\ldots,\frac{p_{n_u}(x)}{q_{n_u}(x)}]^\top$. The dynamical system can be written in the form:
\begin{equation} \label{eq:statespace5}
\begin{aligned}
    \dot{x} &= f(x,u), \\
    q_{k}(x)u_{k} - p_{k}(x) &= 0, \: \forall k = 1, \dots, n_{u},\\
    g_{i}(x,u) &\geq 0, \: \forall i = 1, \dots, q_{1}, \\
    h_{j}(x,u) &= 0, \: \forall j = 1, \dots, q_{2}, \\
    q_{k}(x) &> 0, \: \forall k = 1, \dots, n_{u}, 
\end{aligned}
\end{equation}
where $q_k(x) > 0$ ensures that the controller has no asymptotes and is valid in the state space—this avoids sliding modes; while it is possible to manage sliding modes, this creates difficulties and would require further tests to verify stability.  The following proposition demonstrates how the design problem can be written as an optimisation problem.

\begin{proposition} \label{prop:ss5stability}
Consider System \eqref{eq:statespace5} and suppose there exists a polynomial function $V(x)$ satisfying
\begin{equation} \label{sosopt5}
\begin{aligned}
    &V(x) - \rho(x) \in \Sigma[x], \quad \rho(x) > 0, \\
    &-\frac{\partial V}{\partial x}(x) f(x,u) - \sum_{k=1}^{n_{u}}\lambda_{k}(x,u)(q_{k}(x)u_{k} - p_{k}(x)) -\\
    & \sum_{j=1}^{q_2}t_{j}(x,u)h_{j}(x,u) - \sum_{i=1}^{q_1}s_{i}(x,u)g_{i}(x,u) \in \Sigma[x,u], \\
    &s_{i}(x,u) \in \Sigma[x,u], \: \forall i = 1, \dots, q_{1}, \\
    &t_{j}(x,u) \in \mathbb{R}[x,u], \: \forall j = 1, \dots,
    q_{2}, \\
    &\lambda_{k}(x,u) \in \mathbb{R}[x,u], \: \forall k = 1, \dots, n_{u}, \\
    &p_{k}(x) \in \mathbb{R}[x], \: \forall k = 1, \dots, n_{u} \\
    &q_{k}(x) - \eta_{k}(x) \in \Sigma[x],\quad \eta_{k}(x) >0,  \: \forall k = 1, \dots, n_{u},
\end{aligned}
\end{equation}
then the equilibrium of the system \eqref{eq:statespace5} is stable.
\end{proposition}

\begin{pf}
    This proposition follows immediately from Proposition \ref{prop:ss1stability}, applied to system~\eqref{eq:statespace5}.
\end{pf}
Proposition \ref{prop:ss5stability} provides an SOS framework for nonlinear systems with a decoupled input expression. When the controller $\pi$ is known, the proposed convex formulation can be used to search for a Lyapunov function. Conversely, when a Lyapunov function is given, the same framework enables the synthesis of a stabilising controller. When it comes to joint synthesis, this problem is non-convex because $\lambda(x,u)(q(x)u - p(x))$ contains two terms with decision variables that are multiplied together. However, as we will see below, this formulation offers advantages that can be exploited in an iterative scheme to design rational state feedback controllers for nonlinear systems.

\subsection{Iterative Scheme}
We assume that it is possible to obtain a linear or polynomial controller as an initial controller to start an iterative scheme. This can be achieved through one of many traditional methods, some of which were introduced in Section \ref{sec:rationalintro}. As linear and polynomial functions are a subset of rational functions, once an initial controller 
$u = K(x)$ that stabilises the system is obtained, we can write it as $u = \frac{p^{0}(x)}{q^{0}(x)}$, where $p^{0}(x) = K(x)$, $q^{0}(x) = 1$. Each iteration of the optimisation problem consists of two steps: first, searching for the multiplier $\lambda^{a}(x,u)$; and second, searching for the controller parameters $p^{a}(x)$ and $q^{a}(x)$, where the superscript denotes the $a^{\mathrm{th}}$ iteration. For each step and iteration of the optimisation problem a \emph{new} Lyapunov function is computed, i.e. $V^{a}_1(x)$ and $V^{a}_2(x)$. Here the subscripts 1 and 2 stand for step 1 and step 2 of the $a^{\mathrm{th}}$ iteration respectively.

To formalise the process, consider System \eqref{eq:statespace5} and a known linear or polynomial controller $p^{0}(x)/q^{0}(x)$. A stabilising rational controller can be found at each iteration of a two-step optimisation problem as follows: 

\textbf{Step 1:} with an initial controller, the search for a CLF and multiplier is posed as a feasibility problem:
\begin{equation} \label{opt:iterative1}
\begin{aligned}
& \text{find } V^a_1,\lambda^a_k \qquad \text{such that} \\
    &V^{a}_1(x) - \rho^{a}_1(x) \in \Sigma[x], \quad \rho^{a}_1(x) > 0, \\
    &-\frac{\partial V^{a}_1}{\partial x}(x) f(x,u) - \sum_{k=1}^{n_{u}}\lambda_{k}^{a}(x,u)(q_{k}^{*}(x)u_{k} - p_{k}^{*}(x)) -\\
    & \sum_{j=1}^{q_2}t_{1,j}^{a}(x,u)h_{j}(x,u) - \sum_{i=1}^{q_1}s_{1,i}^{a}(x,u)g_{i}(x,u) \in \Sigma[x,u], \\
    &s_{1,i}^{a}(x,u) \in \Sigma[x,u], \: \forall i = 1, \dots, q_{1}, \\
    &t_{1,j}^{a}(x,u) \in \mathbb{R}[x,u], \: \forall j = 1, \dots,
    q_{2}, \\
    &\lambda_{k}^{a}(x,u) \in \mathbb{R}[x,u], \: \forall k = 1, \dots, n_{u},
\end{aligned}
\end{equation}

\textbf{Step 2:} given the multiplier $\lambda^*_k$, the search for both CLF and a new controller is posed as another feasibility problem:
\begin{equation} \label{opt:iterative2}
\begin{aligned}
& \text{find } V^a_2,p^a_k,q^a_k \qquad \text{such that} \\
    &V^{a}_2(x) - \rho^{a}_2(x) \in \Sigma[x],\quad
    \rho^{a}_2(x) > 0, \\
    &-\frac{\partial V^{a}_2}{\partial x}(x) f(x,u) - \sum_{k=1}^{n_{u}}\lambda_{k}^{*}(x,u)(q_{k}^{a}(x)u_{k} - p_{k}^{a}(x)) -\\
    &\sum_{j=1}^{q_2}t_{2,j}^{a}(x,u)h_{j}(x,u) - \sum_{i=1}^{q_1}s_{2,i}^{a}(x,u)g_{i}(x,u) \in \Sigma[x,u], \\
    &s_{2,i}^{a}(x,u) \in \Sigma[x,u], \: \forall i = 1, \dots, q_{1}, \\
    &t_{2,j}^{a}(x,u) \in \mathbb{R}[x,u], \: \forall j = 1, \dots,
    q_{2}, \\
    &p_{k}^{a}(x) \in \mathbb{R}[x], \: \forall k = 1, \dots, n_{u}, \\
    &q_{k}^{a}(x) - \eta_{k}^{a}(x) \in \Sigma[x], \: \forall k = 1, \dots, n_{u}, \\
    &\eta_{k}^{a}(x) >0,\: \forall k = 1, \dots, n_{u}.
\end{aligned}
\end{equation}
The superscript asterisk represents the solutions from the optimisation problem of the previous step. 
The iteration terminates when the obtained controller solution converges. Alternatively, a maximum number of iterations can be specified to terminate the process early. Furthermore, the degrees of the CLF, multipliers, and the numerator and denominator of the controller in the SOS problem can be increased to obtain a controller that satisfies the desired conditions in, e.g., larger areas around the equilibrium of interest.

\begin{figure}[!ht] 
    \centering  
    \includegraphics[width=1\linewidth]{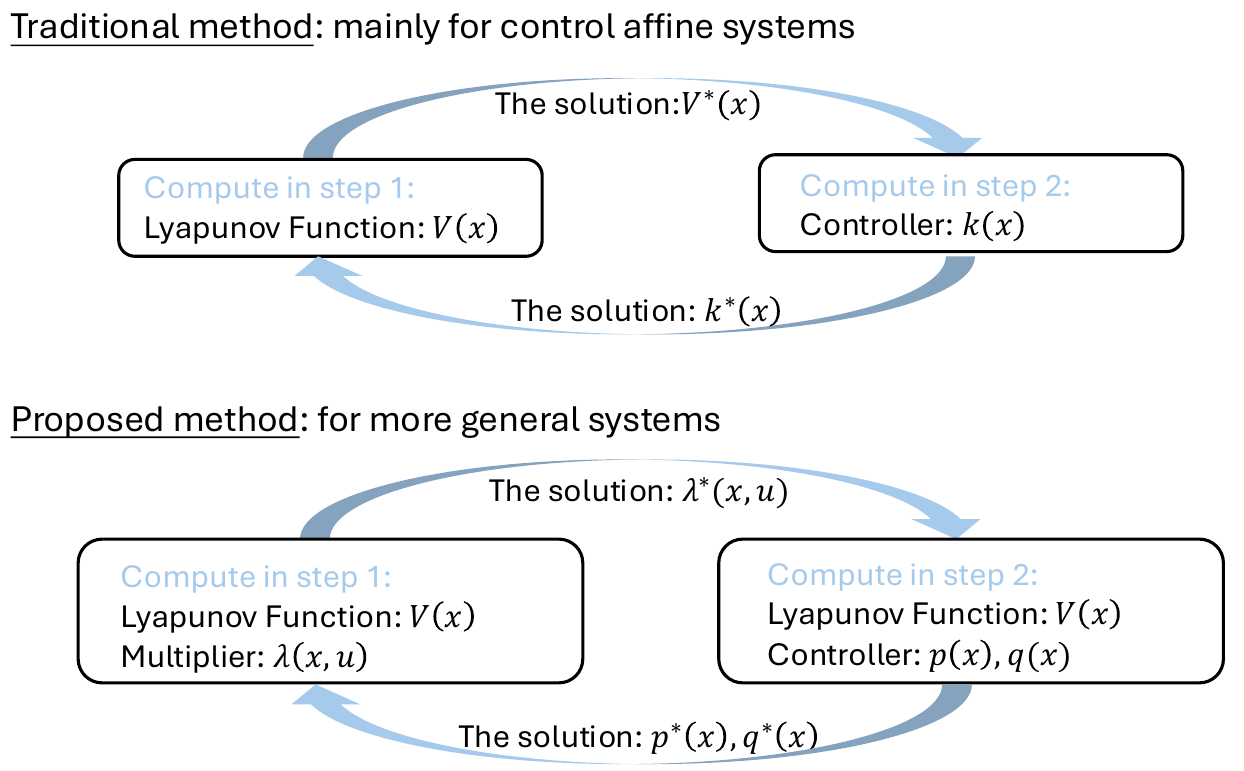}
    \caption{A comparison between traditional iterative approaches and the proposed approach.} \label{methods comparison}
\end{figure}

Fig.~\ref{methods comparison} is a schematic diagram of the solution process of two iterative methods - a traditional iterative approach, and ours. Compared to the proposed method, the traditional iterative method is limited as the controller parameters are not being searched for at the same time as the Lyapunov function. On the contrary, in our method a new Lyapunov is being searched for in both iterative steps. By searching for a new Lyapunov function each time we introduce more flexibility into the optimisation problem and allow a larger variable and parameter space to be searched over. Additionally, by decoupling the constraint and introducing $u$ as a variable in the optimisation problem, we allow further flexibility and expansion of the problem space. 

It is important to note that another significant advantage of our proposed procedure over traditional iterative methods is its applicability to a more general class of nonlinear systems as we noted in Fig.~\ref{methods comparison}. As we illustrate in Section \ref{sec:traditional iterative method}, conventional approaches are typically tailored for control-affine systems, where the control input $u$ appears linearly in the dynamics. For systems with more complex, non-affine dependencies on the control input, these methods can become impractical or even fail. For instance, consider a system with polynomial but non-affine terms in $u$, such as $\dot{x} = f(x) + g_1(x)u + g_2(x)u^2$. In the traditional approach, substituting the rational controller $u=p(x)/q(x)$ and clearing the denominator would introduce complex non-convex cross-terms like $\frac{\partial V}{\partial x}g_1(x)p(x)q(x)$, which involve the product of three different decision variables and are not readily handled by alternating convexification. For dynamics with higher powers of the input, like $\dot{x} = -x + x^2u^3$, this substitution leads to a computational explosion, as the resulting stability conditions would involve high-degree polynomials such as $p(x)^3$ and $q(x)^3$, making the SOS problem difficult to solve. Although the proposed method cannot eliminate the non-convexity, it avoids the aforementioned problems by systematically structuring the problem to maintain tractability. By decoupling the controller as a constraint and introducing a Lagrange multiplier $\lambda(x,u)$, the non-convexity is consistently contained within the bilinear products of the multiplier and the controller polynomials. This structure is crucial because it is always amenable to the two-step iterative convex scheme. By handling the full dynamics $f(x,u)$ directly, our approach isolates the challenging non-convex coupling into a consistent and manageable form, thus making the design problem computationally feasible for a much broader class of nonlinear systems.

Furthermore, as each step in the optimisation problem is convex, our design procedure enables rational controllers to be synthesised even in scenarios where traditional methods may fail to converge or provide tractable solutions. This flexibility allows our method to serve not only as a standalone synthesis approach, but also as a post-processing stage to fine-tune existing controllers. For instance, an initial linear or polynomial controller obtained by another method can be used to warm-start our iteration, which then yields a more expressive rational controller with improved performance. 

\subsection{Extension of Iterative Scheme}
This iterative scheme can be extended in several directions. For instance, we can implement a relevant objective function or add constraints on the system's operating region. Essentially, the proposed method consists of two stages, one is to search for a multiplier and hence to optimise the feasible region of an existing controller and the other is to find the best controller for the objective. Therefore, the objective function of each optimisation step should align with this method. Instead of having both objectives in one optimisation problem, we can split them into the two optimisation problems, allowing for convex design.

As shown in \cite{apap3} we can optimise the controller performance by incorporating the constraint
\begin{equation} \label{eq:decayrate}
    \dot{V}(x) \leq - \gamma V(x),
\end{equation}
where $\gamma$ is a given decay rate that can be increased at each iteration.

A simple way to maximise the region of attraction that the system operates in is by expanding the size of the region at each iteration. For example, we could add a constraint of the form
\begin{equation} \label{eq:ballradius}
    \sum_{i}x_{i}^{2} \leq R,
\end{equation}
where $R$ is the radius of the ball. $R$ can be increased at each iteration. Therefore, increasing $\gamma$ in  \eqref{eq:decayrate} and $R$ in \eqref{eq:ballradius} allows the two objectives to be optimised until one of the optimisation problems becomes infeasible and hence the best performing controller is achieved. The performance objective can be included during the optimisation step to find the rational function parameters, whereas the region objective can be performed at both steps. As the main result of this work, Algorithm \ref{alg:iterative} gives a description of how the proposed iterative scheme is solved to obtain a rational controller. $R^{\text{inc}}$ and $\gamma^{\text{inc}}$ represent the increment of the radius and decay rate after each iteration respectively. $d^{(\cdot)}$ is the degree of polynomial and $d^{(\cdot)}_{\max}$ is the maximum degree that we set. $iter_{\text{max}}$ is the maximum number of iterations that we set. 

\begin{algorithm}
\caption{Convex Iterative Method for Rational Controller Design}
\label{alg:iterative}
\begin{algorithmic}[1]
\State \textbf{Inputs:}
\State \parbox[t]{\linewidth}{Parameters: $R$, $\gamma$, $R^{\text{inc}}$, $\gamma^{\text{inc}}$, $d^{p}$, $d^q$, $d^V$, $d^{\lambda}$, $d^t$, $d^s$, $d^V_{\text{max}}$, $d^{\lambda}_{\text{max}}$, $d^{t}_{\text{max}}$, $d^{s}_{\text{max}}$, $iter_{\text{max}}$}
\State Initialize: $a \gets 0$, $(p^*, q^*) \gets (p^0, q^0)$
\algrule
\State \textbf{Outer Loop:}
\If{($d^{V},d^{\lambda},d^{t},d^{s}) =(d^{V}_{\max},d^{\lambda}_{\max},d^{t}_{\max},d^{s}_{\max})$}
    \State \textbf{STOP} 
\Else
    \If{$a > 0$}
        \State $d^{V} \gets d^{V} + 2$ \textbf{ if } $d^{V} < d_{\max}^{V}$
        \State $d^{\lambda} \gets d^{\lambda} + 2$ \textbf{ if } $d^{\lambda} < d_{\max}^{\lambda}$
        \State $d^{t} \gets d^{t} + 2$ \textbf{ if } $d^{t} < d_{\max}^t$
        \State $d^{s} \gets d^{s} + 2$ \textbf{ if } $d^{s} < d_{\max}^s$
        \State $a \gets 0$
    \EndIf
    \State Go to \textbf{Inner Loop}
\EndIf
\algrule
\State \textbf{Inner Loop:}
\While{$a<iter_{\text{max}}$}
    \State $a \gets a + 1$
    \State \textbf{Step 1.} Search for multiplier and LF:
    \State Solve problem \eqref{opt:iterative1} $+$ \eqref{eq:ballradius} with $(p^*, q^*)$
    \If{the problem is feasible}  
        \State $\lambda^* \gets \lambda^a$, $V^* \gets V^a_1$
    \Else 
        \State Go to \textbf{Outer Loop}
    \EndIf 
    \State \hrulefill
    \State \textbf{Step 2.} Search for controller and CLF:
    \State Solve problem \eqref{opt:iterative2} $+$ \eqref{eq:decayrate} $+$ \eqref{eq:ballradius} with $\lambda^*$
    \If{the problem is feasible}  
        \State $(p^*, q^*) \gets (p^a, q^a)$, $V^* \gets V^a_2$
        \State $R \gets R+R^{\text{inc}}$, $\gamma \gets \gamma+\gamma^{\text{inc}}$
    \Else 
        \State Go to \textbf{Outer Loop}
    \EndIf 
\EndWhile
\algrule
\State \textbf{Outputs:}
\State \Return $(p^*, q^*, V^*)$
\end{algorithmic}
\end{algorithm}

Compared with the standard two-step iterative approach presented in the previous section, Algorithm \ref{alg:iterative} offers a more complete iterative structure by incorporating the optimisable objectives as constraints. We naturally drive the iteration forward by increasing the radius and the decay rate in each iteration. The iteration concludes when one of the optimisation problem becomes infeasible, which indicates that a certain objective can no longer be optimised. Furthermore, we augment this framework with an outer loop: when the problem is infeasible, we can increase the degree of the multipliers and Lyapunov function in the SOS problem. This effectively expands the search space to restore feasibility of the optimisation problem. The final controller is obtained when the degree of all polynomials reaches a pre-set maximum limit (corresponding to computational resource availability) or when the maximum number of inner loop iterations is reached. It is obvious that this iterative method can yield a superior rational controller after each iteration, and if feasibility cannot be restored, the solution returned is identical to the controller in the previous iteration.

The proposed iterative procedure allows the incorporation of additional performance metrics as convex constraints, making the overall design process more versatile. This flexibility also makes it possible to extend the framework with other controller design architectures. One notable example is NNs, which has been motivated by their ability to perform well in numerous machine learning tasks \cite{czhang,bboc}. However, these learnt policies often lack interpretability, may perform poorly when the real environment is different from the learnt environment \cite{vbeh}, and pose challenges for safety verification~\cite{rsut}. In contrast, our rational controller procedure ensures convexity at each optimisation step, and offers better transparency and tractability for analysis. Recent work on rational NNs~\cite{nbou} highlights the potential of combining expressivity with structure, pointing to promising directions for future integration within our procedure.

\section{A Unifying Theoretical Perspective on the Proposed Procedure}
\label{sec:unifying_perspective}

A key strength of our proposed procedure is its generality. This section provides a deeper analysis by demonstrating how this procedure unifies and recovers several well-known rational controller synthesis methods. This analysis reveals that, compared to traditional methods that rely on problem-specific, analytical calculations, our general and systematic iterative numerical approach offers advantages in terms of generality and applicability. We first illustrate how the introduction of the multiplier $\lambda$ enhances the flexibility of the proposed method by explaining its role within the Positivstellensatz framework. Then we show that by imposing specific constraints on this multiplier of our synthesis condition, we can recover the formulations of three distinct classes of controller design methods.

\subsection{The Multiplier as an Algebraic Certificate}
The theoretical foundation of our method is to reframe the synthesis problem using the Positivstellensatz, which yields the SOS stability condition in Proposition \ref{prop:ss5stability}. For simplicity, in this section we consider this stability condition for a single-input system without constraints
\begin{equation}\label{eq:proposed con simplified}
    -\frac{\partial V}{\partial x}(x) f(x,u) - \lambda(x,u)(q(x)u - p(x)) \in \Sigma[x,u].
\end{equation}
This formulation introduces a multiplier polynomial $\lambda(x,u)$, which is analogous to a Lagrange multiplier but is itself an optimisable function within the synthesis problem. The role of this multiplier is to provide an algebraic certificate of stability that is valid for all $(x,u)$.

The key insight is that by treating the coefficients of $\lambda(x,u)$ as decision variables, the optimiser can shape the term $-\lambda(x,u)(q(x)u - p(x))$ to ensure the global non-negativity of the entire SOS expression. While this term is identically zero on the controller-defined set $\mathcal{C} = \{ (x, u) \mid q(x)u - p(x) = 0 \}$, its flexibility off this set is crucial. This approach fundamentally transforms a rigid verification problem, where a fixed expression must be tested for SOS feasibility, into a more flexible search problem. The practical consequence is a significant enlargement of the feasible set for the synthesis problem, as the optimiser actively searches for a multiplier that renders the condition feasible. As will be shown, prior synthesis techniques can be interpreted as implicitly imposing a fixed, highly restrictive structure on this multiplier, thereby constraining the search to a smaller solution space.


\subsection{Connection to Methods based on Analytical Cancellation}
Our method begins with the stability condition \eqref{eq:proposed con simplified} derived from decoupling the controller as a constraint.
To establish a connection with methods based on cancellation, we first consider a control-affine system \eqref{eq:statespace6}. We then constrain the multiplier $\lambda$ to be a function of the state $x$ only, i.e. $\lambda(x)$, and require the certificate of non-negativity to be an SOS polynomial in $x$ only, which we denote as $W(x)$. Under these constraints, the stability condition \eqref{eq:proposed con simplified} becomes an algebraic identity
\begin{equation}
    -\frac{\partial V}{\partial x} f(x) -\frac{\partial V}{\partial x} g(x)u - \lambda(x)(q(x)u - p(x)) - W(x) = 0,
\end{equation}
where $W(x)\in \Sigma[x]$. As this identity must hold for all $u$, and it is linear in $u$, the coefficients of the terms dependent and independent of $u$ must both be zero. This yields two conditions:
\begin{subequations}
    \begin{align}
        -\frac{\partial V}{\partial x}f(x)+\lambda(x)p(x)-W(x) &= 0,\\
        -\frac{\partial V}{\partial x}g(x)-\lambda(x)q(x) &= 0.
    \end{align}
\end{subequations}
Solving these equations for the controller polynomials $p(x)$ and $q(x)$ gives the resulting rational controller
\begin{equation}
    u(x)=\frac{p(x)}{q(x)}=-\frac{\frac{\partial V}{\partial x} f(x)+W(x)}{\frac{\partial V}{\partial x}g(x)},
\end{equation}
which is the same as \eqref{eq:sontag}, the rational controller designed by analytical cancellation inspired by Sontag's universal formula. This demonstrates that the classical analytical result is a special case of our method under a constrained multiplier and certificate structure.

\subsection{Connection to Traditional Iterative Methods}
The Lyapunov stability condition \eqref{eq:traditional sta con} commonly used in traditional iterative design methods can also be recovered. By again considering a control-affine system and making the specific multiplier choice of $\lambda(x)=-\frac{\frac{\partial V}{\partial x}g(x)}{q(x)}$ in the proposed stability condition \eqref{eq:proposed con simplified}, the terms involving $u$ algebraically cancel. This leaves the simplified condition \eqref{eq:traditional sta con}, 
which is exactly the non-convex Lyapunov condition that traditional iterative methods aim to solve. This shows that the traditional approach is equivalent to our method but with the multiplier $\lambda$ implicitly fixed to a rigid, pre-defined form rather than being treated as an optimisable variable.

\subsection{Connection to Completing the Square Methods}
Other methods, such as those by Prajna et al. \cite{prajna2004nonlinear} and Huang et al. \cite{huang2013robust}, achieve a one-shot convexification for systems with a specific structure by employing a strategy analogous to completing the square. To establish the connection, we briefly outline this formulation.

Assume that the system can be expressed in the state-dependent linear-like form $\dot{x} = A(x)Z(x) + B(x)u$ and adopt a matching Lyapunov function of the form $V(x)=Z(x)^\top P^{-1}(\tilde{x})Z(x)$. Define the index set $J = \{j_1, j_2, \dots, j_m\}$ to contain the indices of those rows in $B(x)$ that are equal to zero and define $\tilde{x} = (x_{j_1}, x_{j_2}, \dots, x_{j_m})$ as the subvector of $x$ corresponding to these indices. In this way, their methods eliminate the cross terms between the Lyapunov function and the control input $u$ that may appear in the Lyapunov stability condition, thereby achieving convexification. Through a detailed derivation, the Lyapunov derivative $\dot{V}$ can be shown to have a quadratic structure in the transformed variable $v := P^{-1}(\tilde{x})Z(x)$ and the input $u$
\begin{equation}\label{eq:CS con}
    \dot{V} = v^\top Q(x) v + 2 v^\top S(x) u,
\end{equation}
where the matrices $Q(x)$ and $S(x)$ are functions of the system dynamics and $P(\tilde{x})$. The synthesis goal is to find a controller in the form of $u = K(x)v$ that makes $\dot{V}: = v^\top Q(x) v + 2 v^\top S(x) K(x) v < 0$. This choice results in a convex Lyapunov stability condition of the closed-loop system  
\begin{equation}\label{eq:Prajna Stability}
    -Q(x) - 2 S(x) K(x) \in \Sigma[x].
\end{equation}
It is obvious that our procedure recovers this result by assuming the same Lyapunov structure and selecting the multiplier to be precisely the coefficient of the linear $u$ term in $-\dot{V}$ such that $\lambda = -2v^\top S(x)$.
By choosing this $\lambda$ in the proposed stability condition \eqref{eq:proposed con simplified} and setting the controller $p(x)/q(x)$ to the form derived in \cite{prajna2004nonlinear}, the terms linear in $u$ algebraically cancel, exactly reproducing their final SOS stability condition \eqref{eq:Prajna Stability}. This reveals that their successful but bespoke approach is equivalent to constraining the multiplier $\lambda$ in our method to a highly specific form, which is implicitly defined by their problem structure and change-of-variables approach.

\subsection{Discussion on the Flexibility of the Proposed Method}
Our analysis reveals the ability of our method to recover several existing methods, which achieve convexification through problem-specific strategies, such as analytical cancellation or a change of variables. Within our procedure, each of these tricks corresponds to imposing a specific, rigid constraint on the polynomial $\lambda$.

In contrast, our method treats $\lambda(x, u)$ as a decision variable that is optimised at each step. This allows the stability certificate to be numerically co-designed with the controller, rather than being pre-defined by an analytical technique. By allocating more flexibility to this certificate, our procedure can explore a richer solution space. This offers the potential to find less conservative results compared to methods that rely on a fixed structure, as demonstrated by the numerical examples that follow.


\section{Numerical Examples}
In this section we show how our method can be used to design state-feedback controllers. All simulations are performed on a laptop with Apple M2 chip and 8GB of RAM. The SOS programs are implemented using MATLAB and SOSTOOLS to parse the SOS constraints into an SDP, which is solved using MOSEK~\cite{mosek}.

\subsection{Constrained Inverted Pendulum}
To illustrate the effectiveness of the proposed method, we begin by considering a simple linearised inverted pendulum \cite{hyin1} with a sector constraint, which can be used to handle the nonlinearity. The dynamics and parameters for the system are written as
\begin{equation*}
\begin{aligned}
    &\dot{x}_{1} = x_{2},\quad \dot{x}_{2} = \frac{(mg\ell (x_{1} - x_{3}) - \mu x_{2} + u))}{m\ell^2}, \\
    &x_{3} (\alpha x_{1} - x_{3}) \geq 0,\quad\alpha = 1-\frac{\sin(x_{1,\max})}{x_{1,\max}}.
\end{aligned}
\end{equation*}
Here $x_{1,\max}=R/\sqrt{2}$, where $R$ is from \eqref{eq:ballradius}. We enforce the sector nonlinearity in the Positivstellensatz as described in \cite{mnew5} and also incorporate the decay rate equation \eqref{eq:decayrate} to optimise the performance of the controller such that
\begin{equation*}
    -\frac{\partial V}{\partial x}(x)f(x,u) - \gamma V(x) \in \Sigma[x,u].
\end{equation*}
By incorporating \eqref{eq:decayrate} and \eqref{eq:ballradius}, we increase the radius of the region and the decay rate of Lyapunov function in every iteration. In this case, we choose $R$ and $\gamma$ to start from 1 and 0, respectively. For the $(a+1)^{\text{th}}$ iteration, the update steps are $R^{a+1}=R^a+0.1$ and $\gamma^{a+1}=\gamma^a+0.1$.

\begin{figure}[!ht] 
    \centering  \includegraphics[width=0.75\linewidth]{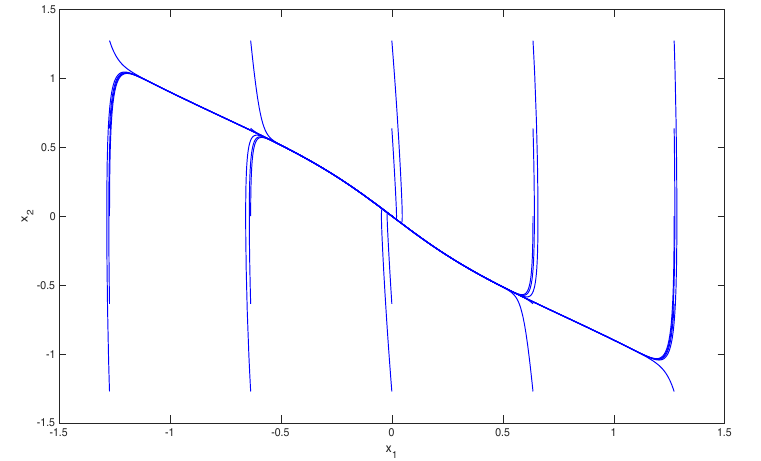}
    \caption{The trajectories under rational controller from the proposed method in the constrained inverted pendulum case.} \label{ConstrainedIP-tra}
\end{figure}

Solving this optimisation problem after 10 iterations provides the rational controller $\frac{p(x_{1},x_{2})}{q(x_{1},x_{2})}$, where $p(x_1,x_2)=
-2.539x_1^3 - 0.66171x_1^2x_2 - 0.52445x_1x_2^2 + 0.00041924x_2^3 - 5.5279x_1 - 1.5831x_2$ and $q(x_1,x_2)=1.1618x_1^2 - 0.69832x_1x_2 + 2.0172x_2^2 + 2.7878$.
The trajectories of the system are shown in Fig.~\ref{ConstrainedIP-tra}. We uniformly select 25 initial points within the considered region - it can be observed that using the controller that we obtained all initial points converge to the equilibrium.

\begin{figure}[!ht] 
    \centering  \includegraphics[width=0.75\linewidth]{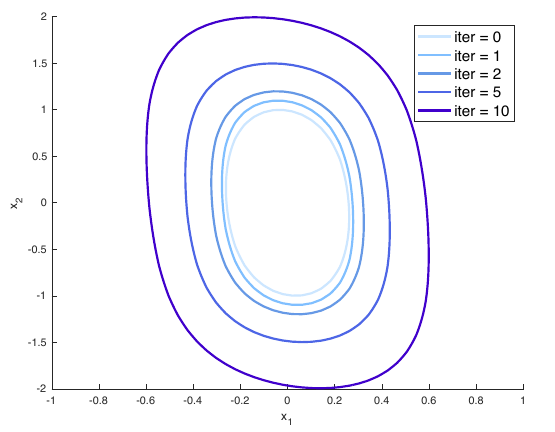}
    \caption{Size of the ROAs at different iteration counts.} \label{ConstrainedIP-ROA}
\end{figure}

The size of the ROAs with rational controllers obtained at different iterations are shown in Fig.~\ref{ConstrainedIP-ROA}. It is clear that because we impose constraints over a larger region with each iteration, a larger ROA is obtained. This illustrates the effectiveness of the proposed rational controller design method. 

\subsection{Rational Function System}
\begin{figure}[!ht] 
    \centering  \includegraphics[width=0.75\linewidth]{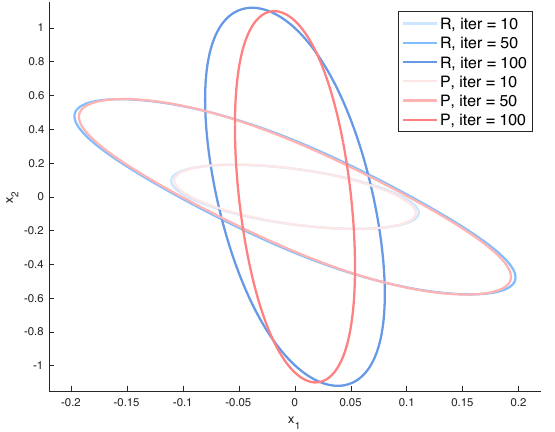}
    \caption{Comparison of ROA size between rational controller and polynomial controller at different iteration counts. R = rational controller, P = polynomial controller.} \label{2D-ROA}
\end{figure}

We now consider a nonlinear system with dynamics represented as rational functions to demonstrate the superiority of rational controllers over polynomial controllers. The system dynamics can be written as:
\begin{equation*}
\begin{aligned}
    \dot{x}_{1} = \frac{1 + x_{1}^{2}}{2}x_{2}, \quad \dot{x}_{2} = \frac{2x_1}{1 + x_{1}^{2}} - x_{2} - \frac{1 - x_{1}^{2}}{1 + x_{1}^{2}}u.
\end{aligned}
\end{equation*}
The proposed controller design method is used to design both a rational controller and a polynomial controller. We start from a region with a radius $R=0.1$ then update radius $R=R+0.01$ and decay rate of Lyapunov function $\gamma=\gamma+0.05$ after each iteration. The ROAs of each controller at iteration 10, 50 and 100 respectively are shown in Fig.~\ref{2D-ROA}. We can see that at each iteration, the ROA of the system under the rational controller which is represented by the blue circles is bigger than that under the polynomial controller. As the number of iterations increases, the difference becomes more pronounced.

On top of the volume of ROAs, we also compare other performance metrics of these two different controllers in terms of the total cost and settling time as Table \ref{2D-performance} shows. We consider the region with radius equal to 0.1 and uniformly select 25 initial points within this region. Here we define the cost as the sum of the states and input squared at all time steps. And we define the settling time as the duration required for all states to reach 1\% of their initial values. With 10 iterations, the performance of the two controllers are close. With 50 iterations, the rational controller is better in both performance metrics. After 100 or more iterations, the total cost under the rational controller is significantly lower with almost the same settling time. 

\begin{table}[ht]
\centering
\caption{Performance comparison of rational (R) and polynomial (P) controllers after different iterations.}
\begin{tabular}{@{}cccccc@{}}
\toprule
Iterations & \multicolumn{2}{c}{Total cost} & \multicolumn{2}{c}{Settling time} \\
\cmidrule(r){2-3} \cmidrule(r){4-5}
           & R & P & R & P \\ \midrule
10         & 1.02                & 1.02                  & 257.05              & 232.55               \\
50         & 1.73                & 1.98                  & 90.65               & 96.61                \\
100        & 31.75               & 99.53                 & 30.79               & 30.25                 \\
150        & 103.65               & 180.20                 & 23.85                & 23.16                 \\
\bottomrule
\label{2D-performance}
\end{tabular}
\end{table}

\begin{figure}[ht]
    \centering
        \begin{subfigure}[b]{0.23\textwidth}
            \centering
            \includegraphics[width=1\textwidth]{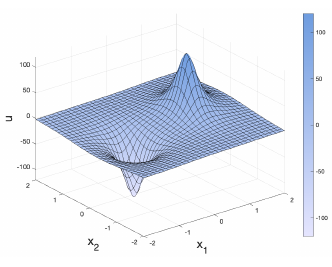}
            \captionsetup{width=\textwidth}
            \caption{Surface of the rational controller, iter = 50}
            \label{2DRationalSurf}
        \end{subfigure}
        \begin{subfigure}[b]{0.23\textwidth}
            \centering
            \includegraphics[width=1\textwidth]{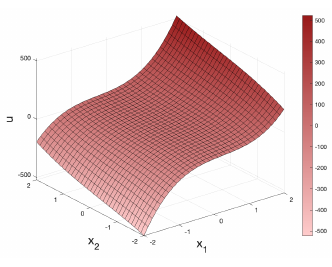}
            \captionsetup{width=\textwidth}
            \caption{Surface of the polynomial controller, iter = 50}
            \label{2DPolynomialSurf}
        \end{subfigure}
    \caption{Controller surface comparison.}
    \label{2DControllerSurf}
\end{figure}

The output of the different controllers with respect to the system states is shown in Fig.~\ref{2DControllerSurf}. Fig.~\ref{2DRationalSurf} and Fig.~\ref{2DPolynomialSurf} illustrate the difference in the output of the rational controller and the polynomial controller after 50 iterations respectively. It is evident that the rational controller exhibits very smooth changes in most cases. One aspect of note is that as the controller approaches its poles the output increases rapidly. In our case, the maximum input value is capped at 100. This should be monitored when applying rational controllers to real tasks. In contrast, the polynomial controller exhibits more drastic variation, especially when far from the equilibrium point, where the required output increases significantly (maximum 500). This also explains the result in Table \ref{2D-performance} - when we increase the number of iterations and consider a larger control range, the total cost required for the system to stabilise under a polynomial controller is significantly greater than that under a rational controller.

\subsection{Nonlinear Polynomial System}
In this section, we consider the polynomial system
\begin{align*}
    \dot{x}_{1} &= -x_{1} + x_{2} - x_{3}, \\\dot{x}_{2} &= -x_{1}(x_{3} + 1) - x_{2}, \\
    \dot{x}_{3} &= -x_{1} + u.
\end{align*}
to demonstrate the superiority of our iterative design method over the traditional iterative method. 

\begin{figure}[ht]
    \centering
    \begin{subfigure}[b]{0.23\textwidth}
        \centering
        \includegraphics[width=\textwidth]{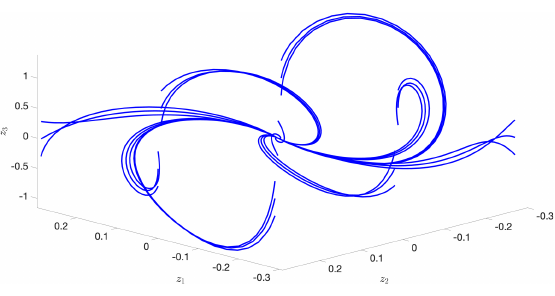}
        \caption{Proposed method}
        \label{3DProposed}
    \end{subfigure}
    \begin{subfigure}[b]{0.23\textwidth}  
        \centering 
    \includegraphics[width=\textwidth]{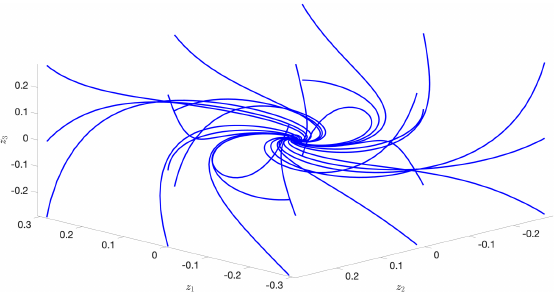}
        \caption{Traditional method}    
        \label{3DTraditional}
    \end{subfigure}
    \caption{Trajectories under controllers obtained from different iterative methods.}
    \label{3DTra}
\end{figure}

We start from a region with radius equal to 0.5 and add 0.1 after each iteration. The proposed iterative method can increase the radius to 6.9 (after 64 iterations) and find a rational controller to stabilise the system within this region, whereas the traditional iterative method cannot find a feasible solution for a region with radius larger than 2.3 (maximum 18 iterations). The trajectories based on the initial states that we select from the initial region under the two rational controllers can be found in Fig.~\ref{3DTra}.

\begin{figure}[ht]
    \centering
    \begin{subfigure}[b]{0.23\textwidth}
        \centering
        \includegraphics[width=1\textwidth]{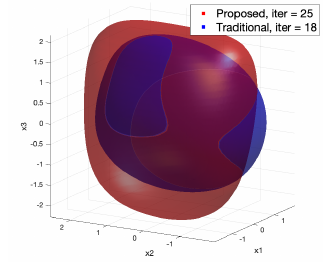}
    \caption{Proposed method: 25 iterations, traditional method: 18 iterations}
    \label{proposed25traditional18}
    \end{subfigure}
    \hfill
    \begin{subfigure}[b]{0.23\textwidth}  
        \centering 
\includegraphics[width=1\textwidth]{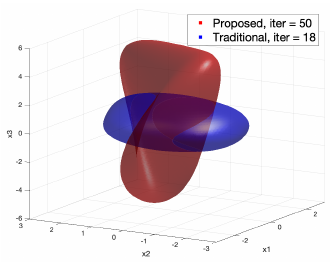}
        \caption{Proposed method: 50 iterations, traditional method: 18 iterations}    
    \label{proposed50traditional18}
    \end{subfigure}
    \caption{Comparison of ROA size under rational controllers obtained from proposed method and traditional method.}
    \label{3DROACompare}
\end{figure}

We also compare the ROA size under the controllers that we obtain from two different methods. We set the iteration number to 18, which is the maximum iteration number of the traditional method and obtain two ROAs of similar shape. The ROA size from the traditional method is slightly bigger than that from the proposed method. However, for the proposed method, we can further iterate and obtain a larger ROA. For example, as Fig.~\ref{proposed25traditional18} shows, when the iteration number is equal to 25, we get a larger ROA compared to the one obtained from the traditional method with 18 iterations. In this case, as we further increase the number of iterations, the ROA is enlarged vertically, as shown in Fig.~\ref{proposed50traditional18}. Although the ROA experiences some horizontal contraction at this point, it significantly enhances the overall range of initial states which can be stabilised. 


We also apply the proposed iterative method to design polynomial controllers by setting $q(x)=1$. Then compare the controller performance with the one obtained from the traditional method. We consider the region with radius equal to 0.5 and uniformly select 25 initial points within this region. The maximum iteration number for both methods is 20. We compare controllers with 10 and 20 iterations respectively. The result is shown in Table \ref{3D-Performance}, where the proposed method demonstrates a better performance for polynomial controller design.

\begin{table}[ht]
\centering
\caption{Performance comparison of polynomial controllers obtained from different methods. P: Proposed, T: Traditional.}
\begin{tabular}{@{}cccccc@{}}
\toprule
Iterations & \multicolumn{2}{c}{Total Cost} & \multicolumn{2}{c}{Settling Time} \\ \cmidrule(lr){2-3} \cmidrule(lr){4-5} 
           & P & T & P & T \\ \midrule
10         & 4.60     & 5.19        & 127.49   & 122.94      \\
20         & 5.61     & 5.96        & 96.68    & 110.98      \\
\bottomrule
\label{3D-Performance}
\end{tabular}
\end{table}

All these results show that we have more flexibility when designing controllers for the system by using the proposed method. The advantages of our method can be explained by the fact that the proposed method does not fix the Lyapunov function. Instead, it solves for the Lyapunov function in both steps 1 and 2 of Algorithm \ref{alg:iterative}, which expands the feasible domain and provides more flexibility for the optimisation problem as we discussed in Section \ref{sec:unifying_perspective}.

\section{Conclusion}
Motivated by rational controllers' strong function approximation capability and inherent smooth saturation property, this paper proposes a new procedure for rational controller design. The key idea is to decouple the controller as a constraint acting on the system. We set up the Lyapunov stability conditions as an SOS problem and solve the non-convex optimisation problem by iteratively solving convex SOS optimisation problems. A stand-out benefit of the proposed iterative method compared to the existing methods is the joint synthesis of the Lyapunov function and the controller -- we search for a new Lyapunov function along with the associated multipliers or controller parameters at each step, instead of alternating between the Lyapunov function and the controller individually. This process not only leads to more flexibility and a richer class of solutions, but as our theoretical analysis reveals, it establishes a uniquely general synthesis framework. We have shown that this inherent certificate flexibility allows our approach to unify a wide range of prior methods. 


Numerical results confirmed the effectiveness of the proposed method for designing rational controllers. Compared with polynomial controllers, rational controllers achieve larger ROA, lower total cost, shorter settling time, and smoother control surface, ensuring that the system can be stabilised by a small input even when deviating significantly from the equilibrium point. Moreover, the proposed convex iterative design method outperforms traditional iterative methods, yielding superior controllers, even when being applied to design polynomial controllers.

This work addresses the rational state-feedback controller design problem by proposing a new iterative procedure. Future directions include incorporating more explicit performance metrics into this design procedure to further enhance controller performance, and combining the advantages of NNs and rational functions and extending the framework to rational NNs \cite{mnew6}.




\bibliographystyle{unsrt}      
\bibliography{autosam}           



\end{document}